\documentclass[draft,preprint,12pt,numbers,sort&compress]{elsarticle}
\usepackage{}
\usepackage{mathrsfs}
\usepackage{amssymb}
\usepackage{amsthm}
\usepackage{mathrsfs}
\usepackage[centertags]{amsmath}
\usepackage{amsfonts}

\usepackage{color}

\usepackage{ctex}
\usepackage{CJK}

\input amssym.def
\input amssym.tex

\date{}

\newtheorem{Theorem}{Theorem}[section]

\newtheorem{Lemma}{Lemma}[section]

\newcommand\R{\mbox{\bf R}}

\newcommand\Z{\mbox{\bf Z}}
\newcommand\T{\mbox{\bf Z}}
\newcommand\z{\mbox{\bf z}}
\newcommand\SR{\mbox{\scriptsize\bf R}}

\newcommand{\definition}{{\lower .5ex
  \hbox{$\>\>\stackrel{\triangle}{=}\>\>$} }}
\newcommand\supp{\mathop{\rm supp}}

\begin{document}

\baselineskip=22pt
\thispagestyle{empty}

\begin{center}
{\Large \bf Convergence problem of  Schr\"odinger equation in
 Fourier-Lebesgue spaces with rough data and random data}\\[1ex]

{ Xiangqian Yan
\footnote{Email:yanxiangqian213@126.com}$^a$, Yajuan Zhao
\footnote{Email:zhaoyj\_91@163.com}$^b$,\qquad Wei Yan\footnote{Email:  011133@htu.edu.cn}$^{a}$}\\[1ex]

{$^a$School of Mathematics and Information Science, Henan
Normal University,}\\
{Xinxiang, Henan 453007,   China}\\[1ex]

{$^b$School of Mathematics, South China  University of  Technology,}\\
{Guangzhou, Guangdong 510640,   China}\\[1ex]

\end{center}
\noindent{\bf Abstract.} In this paper, we consider  the
convergence
 problem of  Schr\"odinger  equation. Firstly,
  we show
 the almost everywhere pointwise convergence
of  Schr\"odinger  equation in Fourier-Lebesgue spaces $\hat{H}^{\frac{1}{p},\frac{p}{2}}(\R)(4\leq p<\infty),$
$\hat{H}^{\frac{3 s_{1}}{p},\frac{2p}{3}}(\R^{2})(s_{1}>\frac{1}{3},3\leq p<\infty),$
$\hat{H}^{\frac{2 s_{1}}{p},p}(\R^{n})(s_{1}>\frac{n}{2(n+1)},2\leq p<\infty,n\geq3)$
   with rough data.
Secondly,  we  show
 that the maximal function
estimate related to  one Schr\"odinger equation  can
 fail with data in  $\hat{H}^{s,\frac{p}{2}}(\R)(s<\frac{1}{p})$.   Finally, we  show the stochastic continuity of
  Schr\"odinger  equation with random data in $\hat{L}^{r}(\R^{n})(2\leq r<\infty)$ almost surely.
The main ingredients are Lemmas 2.4, 2.5, 3.2-3.4.

 \medskip


\medskip
\noindent {\bf Corresponding Author:}Xiangqian Yan

\medskip
\noindent {\bf Email Address:}yanxiangqian213@126.com


\leftskip 0 true cm \rightskip 0 true cm

\baselineskip=20pt

\medskip

{\large\bf 1. Introduction}
\medskip

\setcounter{Theorem}{0} \setcounter{Lemma}{0}\setcounter{Definition}{0}\setcounter{Proposition}{1}

\setcounter{section}{1}

In this paper,  we  investigate  the pointwise convergence  problem of  Schr\"odinger equation
\begin{eqnarray}
&&iu_t+\Delta u=0, \label{1.01}\\
&&u(x,0)=f(x)   \label{1.02}.
\end{eqnarray}
Carleson \cite{Carleson}  showed pointwise
convergence problem of  one dimensional Schr\"odinger equation in
$H^{s}(\R)$ with $s\geq \frac{1}{4}$. Dahlberg and Kenig \cite{DK} showed
that the pointwise convergence of the
Schr\"odinger equation is invalid in $H^{s}(\R^{n})$ with $s<\frac{1}{4}.$
Some authors have studied the  pointwise  convergence  problem of
 Schr\"odinger  equations and  Schr\"odinger type  equations in  higher dimension \cite{Bourgain1995,Bourgain2013,CLV,Cowling,
  DG, DGZ,GS,Lee,LR2015,LR2017,MVV,S,Shao,Tao,TV,Vega,
  Zhang}.
Bourgain \cite{Bourgain2016} presented counterexamples showing that
 when $s<\frac{n}{2(n+1)},n\geq2$,   the pointwise
 convergence problem of $n$ dimensional Schr\"odinger  equation does not hold.
 Du et al. \cite{DGL} investigated  the pointwise convergence problem of
  Schr\"odinger
 equation   in $H^{s}(\R^{2})$ with $s>\frac{1}{3}$.
 Du and Zhang \cite{DZ} showed that the
  pointwise convergence problem of  Schr\"odinger
 equation   is valid  for data in $H^{s}(\R^{n})(s>\frac{n}{2(n+1)},n\geq3).$ Some authors investigated the pointwise
 convergence problem of Schr\"odinger equation on the torus $\mathrm{\T}^{n}$
 \cite{MV,WZ,CLS,EL}.

 Recently, Compaan et al. \cite{CLS} applied randomized initial data method  introduced
  by Lebowitz-Rose-Speer  \cite{LRS}
  and Bourgain  \cite{Bourgain1994,Bourgain1996} and developed by Burq-Tzvetkov \cite{BTL,BTG}
    to study pointwise convergence
   of the Schr\"{o}dinger flow.  For other applications of randomized initial data method,
    we refer the readers to \cite{BOP2015,BOP,CG,CLV,DG,DLM,KMV,KMV2020,M,OOP,OP,ZF2011,ZF2012}.

In this paper, motivated by \cite{KPV1991,Du},  we use the Fourier-Lebesgue space $\hat{H}^{s,r}(\R^{n})$
which is used in \cite{CVV,F,G} to study the pointwise convergence problem of  Schr\"odinger equation. Firstly,
  we show
 the almost everywhere pointwise convergence
of  Schr\"odinger  equation in Fourier-Lebesgue spaces
$\hat{H}^{\frac{1}{p},\frac{p}{2}}(\R)(4\leq p<\infty),$
$\hat{H}^{\frac{3 s_{1}}{p},\frac{2p}{3}}(\R^{2})(s_{1}>\frac{1}{3},3\leq p<\infty),$
$\hat{H}^{\frac{2 s_{1}}{p},p}(\R^{n})(s_{1}>\frac{n}{2(n+1)},2\leq p<\infty,n\geq3)$
with rough data. Secondly, we  present counterexamples showing
 that the maximal function
estimate related to  one Schr\"odinger equation  can
 fail with data in  $\hat{H}^{s,\frac{p}{2}}(\R)(s<\frac{1}{p})$. Finally, motivated by \cite{CLS},
  we show the stochastic continuity of  Schr\"odinger  equation with random  initial data in $\hat{H}^{0,r}(\R)=\hat{L}^{r}(\R^{n})$ almost surely with $2\leq r<\infty$.

We present some notations before stating the main results. $|E|$ denotes by the Lebesgue measure of set $E$.
\begin{eqnarray*}
&&\hat{f}(\xi)=\frac{1}{(2\pi)^{\frac{n}{2}}}\int_{\SR^{n}}e^{-ix\cdot\xi}f(x)dx,
\mathscr{F}_{x}^{-1}f(\xi)=\frac{1}{(2\pi)^{\frac{n}{2}}}\int_{\SR^{n}}e^{ix\cdot\xi}f(x)dx,\\
&&U(t)u_{0}=\frac{1}{(2\pi)^{\frac{n}{2}}}\int_{\SR^{n}}e^{ix\cdot\xi-it|\xi|^{2}}\hat{u}_{0}(\xi)d\xi,
D_{x}^{\alpha}U(t)u_{0}=\frac{1}{(2\pi)^{\frac{n}{2}}}\int_{\SR^{n}} \left|\xi\right|^{\alpha}e^{ix\cdot\xi-it|\xi|^{2}}
\hat{u}_{0}(\xi)d\xi,\\
&&\|f\|_{L_{x}^{q}L_{t}^{p}}=\left(\int_{\SR^{n}}
\left(\int_{\SR}|f(x,t)|^{p}dt\right)^{\frac{q}{p}}dx\right)^{\frac{1}{q}},\frac{1}{r}+\frac{1}{r^{\prime}}=1.
\end{eqnarray*}
$\hat{H}^{s,r}(\R^{n})=\left\{f\in \mathscr{S}^{'}(\R^{n}):\|f \|_{\hat{H}^{s,r}(\SR^{n})}=
\|\langle\xi\rangle^{s}\hat{f}\|_{L_{\xi}^{r^{\prime}}(\SR^{n})}<\infty\right\}$, where
 $\langle\xi\rangle^{s}=(1+|\xi|^{2})^{\frac{s}{2}}$.
\indent Now we introduce the randomization procedure for the initial data, which can be seen
 \cite{BOP2015,BOP,ZF2012}. Let $B^{n}(0,1)=\left\{\xi\in \R^{n}:|\xi|\leq 1\right\}$.
Let  $\psi\in C_{c}^{\infty}(\R^{n})$ be a real-valued,
even, non-negative bump function with $\supp\psi\subset B^{n}(0,1)$ as well as for all $\xi\in\R^{n}$, such that
\begin{eqnarray}
&&\sum\limits_{k\in \z^{n}}\psi(\xi-k)=1,\label{1.03}
\end{eqnarray}
which is known as  Wiener decomposition of  the frequency space. For every $k\in \Z^{n}$,
 we define the function $\psi(D-k)f:\R^{n}\rightarrow\mathbb{C}$ by
\begin{eqnarray}
(\psi(D-k)f)(x)=\mathscr{F}^{-1}\big(\psi(\xi-k)\mathscr{F}f\big)(x),x\in \R^{n}\label{1.04}.
\end{eqnarray}
In fact the above projections satisfy a unit-scale Bernstein inequality, namely that for
 $p_{1}$ and $p_{2}$ satisfying the conditions $2 \leq p_1 \leq p_2 \leq\infty$, from Lemma 2.1 of \cite{LMCPDE}, we know that there
 exists a $C\equiv C(p_1, p_2)>0$ such that for any $f\in L_{x}^2(\R^{n})$ and for any $k\in \Z^{n}$
\begin{eqnarray}
\left\|\psi(D-k)f\right\|_{L_{x}^{p_2}(\SR^{n})}\leq C
\left\|\psi(D-k)f\right\|_{L_{x}^{p_1}(\SR^{n})}\leq C\left\|\psi(D-k)f\right\|_{L_{x}^{2}(\SR^{n})}.\label{1.05}
\end{eqnarray}
Let $\{g_k\}_{k\in \z
^n}$ be a sequence of independent, zero-mean, real-valued  random variables
 on a probability space $(\Omega,\mathcal{A}, \mathbb{P})$,
  endowed with probability
distributions $\mu_n(n\geq 1)$, respectively.
 Assume that $\mu_{n}$ satisfy the property
\begin{eqnarray}
&&\exists c>0,\forall \gamma\in \R, \forall n\geq1, \Big|\int_{-\infty}^{+\infty}e^{\gamma x}d\mu_n(x)\Big|\leq e^{c\gamma^2}.\label{1.06}
\end{eqnarray}
Obviously, (\ref{1.06})  is satisfied by standard Gaussian random examples.
Then, for $f\in \hat{H}^{0,r}(\R^{n})$, we define its randomization by
\begin{eqnarray}
&&\hat{f}^{\omega}:=\sum\limits_{k\in\z^{n}}g_{k}(\omega)\psi(\xi-k)\hat{f}.\label{1.07}
\end{eqnarray}
We define
$
\|f\|_{L_{\omega}^{p}(\Omega)}=\left[\int_{\Omega}|f(\omega)|^{p}d\mathbb{P}(\omega)\right]^{\frac{1}{p}}.
$

The main results are as follows:
\begin{Theorem} \label{Theorem1}
Let  $f\in  \hat{H}^{\frac{1}{p},\frac{p}{2}}(\R)$ with $4\leq p<\infty.$
 Then, we have
\begin{eqnarray}
\lim\limits_{t\longrightarrow 0}U(t)f(x)=f(x)\label{1.09}
\end{eqnarray}
 almost  everywhere.

\end{Theorem}

\noindent{\bf Remark 1.} When $p=4,$ from Theorem 1.1, we obtain the same
result as its of Carleson \cite{Carleson}.
Thus,    we extend the result of \cite{Carleson}.
\begin{Theorem}\label{Theorem2}
Let  $f\in \hat{H}^{\frac{3 s_{1}}{p},\frac{2p}{3}}(\R^{2})$ with $s_{1}>\frac{1}{3},3\leq p<\infty.$
 Then, we have
\begin{eqnarray}
\lim\limits_{t\longrightarrow 0}U(t)f(x)=f(x)\label{1.010}
\end{eqnarray}
 almost  everywhere.
\end{Theorem}
\noindent{\bf Remark 2.} When $p=3,$ from Theorem 1.2, we obtain the same
result as its of  Du et al. \cite{DGL}.
Thus,    we extend the result of \cite{DGL}.
\begin{Theorem}\label{Theorem3}
Let  $f\in \hat{H}^{\frac{2 s_{1}}{p},p}(\R^{n})$ with $s_{1}>\frac{n}{2(n+1)},2\leq p<\infty,n\geq3.$
 Then, we have
\begin{eqnarray}
\lim\limits_{t\longrightarrow 0}U(t)f(x)=f(x)\label{1.011}
\end{eqnarray}
 almost  everywhere.
\end{Theorem}
\noindent{\bf Remark 3.} When $p=2,$ from Theorem 1.3, we obtain the same
result as its of Du and Zhang \cite{DZ}.
Thus,    we extend the result of \cite{DZ}.

\begin{Theorem} \label{Theorem4}
The maximal function  inequality
\begin{eqnarray}
\|U(t)f\|_{L_{x}^{p}L_{t}^{\infty}}
\leq C\|f\|_{\hat{H}^{s,\frac{p}{2}}(\SR)}\label{1.012}
\end{eqnarray}
does not hold
if $s<\frac{1}{p}.$
\end{Theorem}

\begin{Theorem}\label{Theorem5}
For $f\in \hat{L}^{r}(\R^{n}) ,r\geq2$,  we denote by $f^{\omega}$ the randomization of $f$ as
defined in (\ref{1.07}). Then,  $\forall \alpha>0$,  we have
$
\lim\limits_{t\longrightarrow0}\mathbb{P}\left(\omega\in \Omega:\left|U(t)f^{\omega}(x)-f^{\omega}(x)\right|>\alpha\right)=0
$
uniformly with respect to $x$.
More exactly,  $\forall\alpha>0$, $\forall\epsilon>0$
 such that $Ce\epsilon\left(\ln\frac{C_{1}}{\epsilon}\right)^{\frac{1}{2}}<\alpha$  and
when $|t|<\epsilon$, we have
$\mathbb{P}\left(\left\{\omega\in\Omega:|U(t)f^{\omega}-f^{\omega}|>\alpha\right\}\right)\leq3\epsilon.$
\end{Theorem}

The rest of the paper is arranged as follows. In Section 2,  we give some
preliminaries. In Section 3, we present  probabilistic estimates of some random series. In Section 4,
 we prove
 Theorem 1.1. In Section 5, we prove
 Theorems 1.2, 1.3.In Section 6, we prove
 Theorem 1.4.In Section 7, we prove
 Theorem 1.5.

\bigskip

\setcounter{section}{2}

\noindent{\large\bf 2. Preliminaries }

\setcounter{equation}{0}

\setcounter{Theorem}{0}

\setcounter{Lemma}{0}

\setcounter{section}{2}
In this section, we present some preliminaries.

\begin{Lemma}\label{lem2.1}
 Let $f\in L^{2}(\R)$.
 Then, we have
\begin{eqnarray}
&&\left\|U(t)f\right\|_{L_{x}^{4}L_{t}^{\infty}}
\leq \|D_{x}^{\frac{1}{4}}f\|_{L^{2}(\SR)}=\||\xi|^{\frac{1}{4}}\hat{f}\|_{L^{2}(\SR)}.\label{2.01}
\end{eqnarray}
\end{Lemma}

For  the proof of Lemma 2.1, we refer the readers to \cite{KPV1991}.
\begin{Lemma}\label{lem2.2}
Let $f\in H^{s_{1}}(\R^{2})$, $s_{1}>\frac{1}{3}$.Then, we have
\begin{eqnarray}
&&\left\|U(t)f\right\|_{L_{x}^{3}(B(0,1))L_{t}^{\infty}}\leq C\left\|f\right\|_{H^{s}(\SR^{2})}.\label{2.02}
\end{eqnarray}
\end{Lemma}

For the proof of  Lemma 2.2, we refer  the readers  to \cite{DGL}.
\begin{Lemma}\label{lem2.3}
Let $f\in H^{s_{2}}(\R^{n}),n\geq3$, $s_{2}>\frac{n}{2n+1}$. Then, we have
\begin{eqnarray}
&&\left\|U(t)f\right\|_{L_{x}^{2}(B(0,1))L_{t}^{\infty}}\leq C\left\|f\right\|_{H^{s}(\SR^{n})}.\label{2.03}
\end{eqnarray}
\end{Lemma}

For the proof of  Lemma 2.3, we refer  the reader to \cite{DZ}.
\begin{Lemma}\label{lem2.4}(Maximal function estimate related to the Fourier-Lebesgue space)
\begin{eqnarray}
&&\left\|U(t)f\right\|_{L_{x}^{p}(\SR)L_{t}^{\infty}}
\leq C\|f\|_{\hat{H}^{\frac{1}{p},\frac{p}{2}}},4\leq p<\infty\label{2.04},\\
&&\left\|U(t)f\right\|_{L_{x}^{p}(B^{2}(0,1))L_{t}^{\infty}}
\leq C\left\|f\right\|_{\hat{H}^{\frac{3s_{1}}{p},\frac{2p}{3}}},3\leq p<\infty,s_{1}>\frac{1}{3},\label{2.05}\\
&&\left\|U(t)f\right\|_{L_{x}^{p}(B^{n}(0,1))L_{t}^{\infty}}
\leq C\left\|f\right\|_{\hat{H}^{\frac{2s_{2}}{p},p}},2\leq p<\infty,s_{2}>\frac{n}{2(n+1)},n\geq3\label{2.06}.
\end{eqnarray}
\end{Lemma}
\noindent{\bf Proof.}Obviously, we have
\begin{eqnarray}
&&\left\|U(t)f\right\|_{L_{x}^{\infty}(B(0,1))L_{t}^{\infty}}\leq C\left\|U(t)f\right\|_{L_{x}^{\infty}L_{t}^{\infty}}
\leq C\|\hat {f}\|_{L^{1}}.\label{2.07}
\end{eqnarray}
Interpolating (\ref{2.07}) with (\ref{2.01})-(\ref{2.03}) respectively, yields
that (\ref{2.04})-(\ref{2.06}) are valid.

This completes the proof of Lemma 2.4.

\begin{Lemma}\label{lem2.5}(Density Theorem in $\hat{H}^{s,r}(\R^{n})$)
  Let $f \in \hat{H}^{s,r}(\R^{n}), s\in \R,2\leq r<\infty$.
  Then, $\forall \epsilon>0$, there exists a rapidly decreasing function $g$
   and a function $h$ with $\|h\|_{\hat{H}^{s,r}(\SR^{n})}<\epsilon$
  such that $
f=g+h.
$
\end{Lemma}
\noindent{\bf Proof.}From $f\in \hat{H}^{s,r}(\R^{n}),$ we have that $\langle \xi\rangle^{s} \hat {f}\in L^{r^{\prime}}(\R^{n})$, according to
the density theorem in $L^{r^{\prime}}(\R^{n})$, we know that there exists rapidly decreasing function $g_{1}$ and a function $h_{1}$
 such that $\langle \xi\rangle^{s}\hat{f}=g_{1}+h_{1}$ with $\|h_{1}\|_{L^{r^{\prime}}(\SR^{n})}<\epsilon$.
Thus, we have
 \begin{eqnarray}
  f=\mathscr{F}_{x}^{-1}\left(\langle \xi\rangle^{-s}g_{1}\right)+\mathscr{F}_{x}^{-1}\left(\langle \xi\rangle^{-s}h_{1}\right).\label{2.09}
  \end{eqnarray}
Let $g=\mathscr{F}_{x}^{-1}\left(\langle \xi\rangle^{-s}g_{1}\right),h=\mathscr{F}_{x}^{-1}\left(\langle \xi\rangle^{-s}h_{1}\right).$
Since $g_{1}$ is a decreasing rapidly function, thus, $g$ is a decreasing rapidly  function. Obviously,
$
\|h\|_{\hat{H}^{s,r}(\SR^{n})}=\|h_{1}\|_{L^{r^{\prime}}(\SR^{n})}.
$
Thus, we have $f=g+h$. This completes the proof of Lemma 2.5.

\begin{Lemma}\label{lem2.6}
  Let $f$ be a rapidly decreasing function. Then,
  we have
\begin{eqnarray}
\left|U(t)f-f\right|\leq C|t|.\label{2.010}
\end{eqnarray}
\end{Lemma}
For the proof of Lemma 2.6, we refer the readers to pages 14-15 of  Lemma 2.3 in \cite{Du}.

\begin{Lemma}\label{lem2.07}
Let $f\in L^{1}(\R^{n})$, $1\leq q\leq\infty$. Then, we have
\begin{eqnarray}
&&\left[\sum\limits_{N}\left\|P_{N}f\right\|_{L^{1}(\SR^{n})}^{q}\right]^{\frac{1}{q}}\leq \int_{\SR^{n}}\left(\sum\limits_{N}|P_{N}f|^{q}\right)^{\frac{1}{q}}dx.\label{2.012}
\end{eqnarray}
\end{Lemma}
\noindent{\bf Proof}
When $q=\infty$, (\ref{2.012}) is obvious. Then we consider $1\leq q<\infty$.\\
Obviously,
\begin{eqnarray}
&&\left[\sum\limits_{N}\left\|P_{N}f\right\|_{L^{1}}^{q}\right]^{\frac{1}{q}}
=\left[\sum\limits_{N}\left(\int_{\SR^{n}}|P_{N}f|dx\right)^{q}\right]^{\frac{1}{q}}.\label{2.013}
\end{eqnarray}
By duality,  (\ref{2.013}) can be rewritten as
\begin{eqnarray*}
&&\left[\sup\limits_{\|a_{N}\|_{l^{q^{\prime}}}=1}\sum\limits_{N}\left[\int_{\SR^{n}}|P_{N}f|dx\right]a_{N}\right]
=\left[\sup\limits_{\|a_{N}\|_{l^{q^{\prime}}}=1}\int_{\SR^{n}}\sum\limits_{N}|P_{N}f|a_{N}dx\right]\\
&&\leq\left[\sup\limits_{\|a_{N}\|_{l^{q^{\prime}}}=1}\int_{\SR^{n}}\left[\sum\limits_{N}|P_{N}f|^{q}\right]^{\frac{1}{q}}
\|a_{N}\|_{l^{q^{\prime}}}dx\right]\leq\int_{\SR^{n}}\left(\sum\limits_{N}|P_{N}f|^{q}\right)^{\frac{1}{q}}dx.
\end{eqnarray*}
Here  $\frac{1}{q}+\frac{1}{q^{\prime}}=1$.
This completes the proof of Lemma 2.7.
\begin{Lemma}\label{lem2.08}
Let $f\in \hat{L}^{r}(\R^{n})$, $2\leq r<\infty$. Then, we have
\begin{eqnarray}
&&\left[\sum\limits_{k\in\z^{n}}|\psi(D-k)f|^{2}\right]^{\frac{1}{2}}\leq \left\|f\right\|_{\hat{L}^{r}(\SR^{n})}.\label{2.014}
\end{eqnarray}
\end{Lemma}
\noindent{\bf Proof.}
Obviously (\ref{2.014}) is equivalent to
\begin{eqnarray}
&&\sum\limits_{k\in\z^{n}}|\psi(D-k)f|^{2}\leq\left\|f\right\|_{\hat{L}^{r}(\SR^{n})}^{2}.\label{2.015}
\end{eqnarray}
Using the H\"{o}lder inequality, Lemma 2.7 and $\supp\psi\subset B^{n}(0,1)$, $\frac{1}{r}+\frac{1}{r^{\prime}}=1,r\geq2$, $l^{1}\hookrightarrow l^{2},$ since
\begin{eqnarray}
&&|\mathscr{F}_{x}f(\xi)|=\sum\limits_{k\in\z^{n}}\psi(\xi-k)|\mathscr{F}_{x}f(\xi)|,\label{2.016}
\end{eqnarray}
we obtain
\begin{eqnarray}
&&\sum\limits_{k\in\z^{n}}|\psi(D-k)f|^{2}=\frac{1}{(2\pi)^{\frac{n}{2}}}\sum\limits_{k\in \z^{n}}\left|\int_{\SR^{n}}e^{ix\cdot\xi}\psi(\xi-k)\mathscr{F}_{x}f(\xi)d\xi\right|^{2}\nonumber\\
&&=\frac{1}{(2\pi)^{\frac{n}{2}}}\sum\limits_{k\in \z^{n}}\left|\int_{|\xi-k|\leq1}e^{ix\cdot\xi}\psi(\xi-k)\mathscr{F}_{x}f(\xi)d\xi\right|^{2}\nonumber\\
&&\leq C\sum\limits_{k\in\z^{n}}\left[
\left[\int_{|\xi-k|\leq1}\left|\psi(\xi-k)\mathscr{F}_{x}f(\xi)\right|^{r^{\prime}}d\xi\right]^{\frac{2}{r^{\prime}}}
\left[\int_{|\xi-k|\leq1}d\xi\right]^{\frac{2}{r}}\right]\nonumber\\
&&\leq C\sum\limits_{k\in\z^{n}}\left[\int_{|\xi-k|\leq1}\left|\psi(\xi-k)\mathscr{F}_{x}f(\xi)\right|^{r^{\prime}}d\xi
\right]^{\frac{2}{r^{\prime}}}\nonumber\\
&&\leq C\sum\limits_{k\in\z^{n}}\left\|\psi(\xi-k)\mathscr{F}_{x}f(\xi)\right\|_{L^{r^{\prime}}}^{2}
\leq C\left\|\|\psi(\xi-k)\mathscr{F}_{x}f(\xi)\|_{l^{2}}\right\|_{L^{r^{\prime}}}^{2}\nonumber\\
&&\leq C\left\|\sum\limits_{k\in\z^{n}}\psi(\xi-k)|\mathscr{F}_{x}f(\xi)|\right\|_{L^{r^{\prime}}}^{2}\leq C\left\|f\right\|_{\hat{L}^{r}(\SR^{n})}^{2}.\label{2.017}
\end{eqnarray}

This completes the proof of Lemma 2.8.

\begin{Lemma}\label{lem2.9}
Let $f\in \hat{L}^{r}(\R^{n})$, $2\leq r<\infty$. Then, we have
\begin{eqnarray}
&&\left[\sum\limits_{k\in\z^{n}}|\psi(D-k)U(t)f|^{2}\right]^{\frac{1}{2}}\leq C\left\|f\right\|_{\hat{L}^{r}(\SR^{n})}.\label{2.018}
\end{eqnarray}
\end{Lemma}

Lemma 2.9 can be proved similarly to Lemma 2.8.

\bigskip

\noindent{\large\bf 3. Probabilistic estimates of some random series}

\setcounter{equation}{0}

\setcounter{Theorem}{0}

\setcounter{Lemma}{0}
\setcounter{section}{3}
\begin{Lemma}\label{lem3.1}
We assume that (\ref{1.06}) holds. Then, there exists $C>0$ such that
\begin{eqnarray}
\left\|\sum_{k\in\z^n}g_k(\omega)c_k\right\|_{L_{\omega}^p(\Omega)}
\leq C\sqrt{p}\left\|c_k\right\|_{l^{2}(\z^n)},\label{3.01}
\end{eqnarray}
for all $p\geq2 $
and $\{c_k\}\in l^{2}(\Z^n)$.
\end{Lemma}

For the proof of  Lemma 3.1,  we refer the reader  to Lemma 3.1 of \cite{BTL}.
\begin{Lemma}\label{lem3.2}
We assume that $g$ is a rapidly decreasing  function and $g^{\omega}$ the randomization of $g$ be
defined as in (\ref{1.07}).
 Then, $\forall \alpha>0$,  there exist $C>0,C_{1}>0$ such that
\begin{eqnarray}
\mathbb{P}\left(\Omega_{1}^{c}\right)\leq C_{1}e^{-\left(\frac{\alpha}{C|t|e}\right)^{2}},\label{3.02}
\end{eqnarray}
where
$
\Omega_{1}^{c}=\left\{\omega\in \Omega:
\left|U(t)g^{\omega}-g^{\omega}\right|>\alpha\right\}.
$
\end{Lemma}
\noindent{\bf Proof.} By  using Lemma 3.1 and  the H\"older inequality with respect to $\xi$,
 since $g$  is  a  rapidly  decreasing  function and $|e^{it|\xi|^{2}}-1|\leq |t\xi^{2}|,$ by using Lemma 2.7 and 
  $l^{1}\hookrightarrow l^{2}$ as well as (\ref{2.016}),
we have
\begin{eqnarray}
&&\hspace{-2cm}\left\|U(t)g^{\omega}-g^{\omega}\right\|_{L_{\omega}^{p}(\Omega)}\leq C
\sqrt{p}\left[\sum_{k}\left|\int_{\SR^{n}}(e^{it|\xi|^{2}}-1)
e^{ix \xi}\psi(\xi-k)\mathscr{F}g(\xi)d\xi\right|^{2}\right]^{\frac{1}{2}}\nonumber\\
&&\hspace{-2cm}\leq C|t|\sqrt{p}\left[\sum_{k}\left|\int_{|\xi-k|\leq1}|\xi|^{2}\left|\psi(\xi-k)
\mathscr{F}g(\xi)\right|d\xi\right|\right]^{\frac{1}{2}}\nonumber\\
&&\hspace{-2cm}\leq C|t|\sqrt{p}\left[\sum_{k}\left[\int_{|\xi-k|\leq1}|\xi|^{2r^{\prime}}
\left|\psi(\xi-k)\mathscr{F}g(\xi)\right|^{r^{\prime}}d\xi\right]^{\frac{2}{r^{\prime}}}
\left[\int_{|\xi-k|\leq1}d\xi\right]^{\frac{2}{r}}\right]^{\frac{1}{2}}\nonumber\\
&&\hspace{-2cm}\leq C|t|\sqrt{p}\left[\sum_{k}\left[\int_{|\xi-k|\leq1}|\xi|^{2r^{\prime}}
\left|\psi(\xi-k)\mathscr{F}g(\xi)\right|^{r^{\prime}}d\xi\right]^{\frac{2}{r^{\prime}}}\right]^{\frac{1}{2}}\nonumber\\
&&\hspace{-2cm}\leq  C|t|\sqrt{p}\left\|\||\xi|^{2}\psi(\xi-k)\mathscr{F}_{x}g(\xi)\|_{l_{k}^{2}}\right\|_{L^{r^{\prime}}}\leq C|t|\sqrt{p}\left\|\sum\limits_{k\in\z^{n}}|\xi|^{2}\psi(\xi-k)|\mathscr{F}_{x}g(\xi)|\right\|_{L^{r^{\prime}}}\nonumber\\
&&\hspace{-2cm}\leq C|t|\sqrt{p}\left\|g\right\|_{\hat{H}^{2,r}(\SR^{n})}\leq C|t|\sqrt{p}.\label{3.03}
\end{eqnarray}
Thus, by using Chebyshev inequality, from (\ref{3.03}), we have
\begin{eqnarray}
&&\mathbb{P}\left(\Omega_{1}^{c}\right)\leq \int_{\Omega_{1}^{c}}\left[\frac{\left|U(t)g^{\omega}-g^{\omega}
\right|}{\alpha}\right]^{p}d\mathbb{P}(\omega)
\leq \left(\frac{C\sqrt{p}|t|}{\alpha}\right)^{p}\label{3.04}.
\end{eqnarray}
Take
\begin{eqnarray}
p=\left(\frac{\alpha}{Ce|t|}\right)^{2}\label{3.05}.
\end{eqnarray}
If $p\geq2$, from (\ref{3.04}),
then we have
\begin{eqnarray}
&&\mathbb{P}\left(\Omega_{1}^{c}\right)\leq e^{-p}=e^{-\left(\frac{\alpha}{Ce|t|}\right)^{2}}\label{3.06}.
\end{eqnarray}
If $p\leq2$, from (\ref{3.04}),
 we have
\begin{eqnarray}
\mathbb{P}(\Omega_{1}^{c})\leq 1\leq e^{2}e^{-2}\leq C_{1}e^{-\left(\frac{\alpha}
{Ce|t|}\right)^{2}}.\label{3.07}
\end{eqnarray}
Here $C_{1}=e^{2}.$
Thus, from (\ref{3.06}), (\ref{3.07}), we have
\begin{eqnarray}
\mathbb{P}(\Omega_{1}^{c})\leq  C_{1} e^{
\left[-\left(\frac{\alpha}
{Ce|t|}\right)^{2}\right]}.\label{3.08}
\end{eqnarray}

This completes the proof of Lemma 3.2.

\begin{Lemma}\label{lem3.3}
For $h\in \hat{L}^{r}(\R^{n}) ,r\geq2$, $h^{\omega}$ the randomization of $h$ be
defined as in (\ref{1.07}). Then, for $\forall \alpha>0$, there exist $C>0,C_{1}>0$ such that
\begin{eqnarray}
&&\mathbb{P}(\Omega_{2}^{c})\leq C_{1}e^{-\left(\frac{\alpha}{Ce\|h\|_{\hat{L}^{r}}}\right)^{2}},\label{3.09}
\end{eqnarray}
where
$
\Omega_{2}^{c}=\left\{\omega\in\Omega:|U(t)h^{\omega}|>\alpha\right\}.
$
\end{Lemma}
\noindent{\bf Proof.}
Using Lemma 3.1 and Lemma 2.9, we have
\begin{eqnarray}
&&\left\|U(t)h^{\omega}\right\|_{L_{\omega}^{p}(\Omega)}=
\left\|\sum\limits_{k\in \z^{n}}g_{k}(\omega)U(t)\psi(D-k)h\right\|_{L_{\omega}^{p}(\Omega)}\nonumber\\
&&\leq C\sqrt{p}\left[\sum\limits_{k\in\z^{n}}|U(t)\psi(D-k)h|^{2}\right]^{\frac{1}{2}}\leq C\sqrt{p}\left\|h\right\|_{\hat{L}^{r}}.\label{3.010}
\end{eqnarray}
By using a proof similar to  (\ref{3.04})-(\ref{3.08}) and (\ref{3.010}),  we obtain that (\ref{3.09}) is valid.

\begin{Lemma}\label{3.4}
For $h\in \hat{L}^{r}(\R^{n}) ,2\leq r<\infty$, $h^{\omega}$ the randomization of $h$ be
defined as in (\ref{1.07}). Then, for $\forall \alpha>0$, there exist $C>0,C_{1}>0$ such that
\begin{eqnarray}
&&\mathbb{P}(\Omega_{3}^{c})\leq C_{1}e^{-\left(\frac{\alpha}{Ce\|h\|_{\hat{L}^{r}}}\right)^{2}},\label{3.011}
\end{eqnarray}
where
$
\Omega_{3}^{c}=\left\{\omega\in\Omega:|h^{\omega}|>\alpha\right\}.
$
\end{Lemma}
\noindent{\bf Proof.}
Using Lemma 3.1 and Lemma 2.8, we have
\begin{eqnarray}
&&\left\|h^{\omega}\right\|_{L_{\omega}^{p}(\Omega)}=\left\|\sum\limits_{k\in \z^{n}}g_{k}(\omega)\psi(D-k)h\right\|_{L_{\omega}^{p}(\Omega)}\nonumber\\
&&\leq C\sqrt{p}\left[\sum\limits_{k\in\z^{n}}|\psi(D-k)h|^{2}\right]^{\frac{1}{2}}\leq C\sqrt{p}\left\|h\right\|_{\hat{L}^{r}}.\label{3.012}
\end{eqnarray}
By using a proof similar to  (\ref{3.04})-(\ref{3.08}) and (\ref{3.012}),  we obtain that (\ref{3.011}) is valid.

This completes the proof of Lemma 3.4.

\medskip

\noindent{\large\bf 4. Proof of Theorem 1.1}
\setcounter{equation}{0}
\setcounter{Theorem}{0}

\setcounter{Lemma}{0}

\setcounter{section}{4}
In this section, we apply (\ref{2.04}),  Lemmas 2.5-2.6 to establish Theorem 1.1.

\noindent{\bf  Proof of Theorem 1.1.}
If $f$ is rapidly decreasing function,
from Lemma 2.6, we have
\begin{eqnarray}
\left|U(t)f-f\right|\leq C|t|.\label{4.01}
\end{eqnarray}
From (\ref{4.01}), we know that Theorem 1.1  is valid.

When $f\in \hat{H}^{\frac{1}{p},\frac{p}{2}}(\R)$,
by using Lemma 2.5, there exists a rapidly decreasing function
 $g$ such that $f=g+h$, where $\|h\|_{\hat {H}^{\frac{1}{p},\frac{p}{2}}(\SR)}<\epsilon.$
Thus, we have
\begin{eqnarray}
\lim\limits_{t\longrightarrow0}\left|U(t)f-f\right|\leq
\lim\limits_{t\longrightarrow0}\left|U(t)g-g\right|
+\lim\limits_{t\longrightarrow0}\left|U(t)h-h\right|
.\label{4.02}
\end{eqnarray}
We  define
\begin{eqnarray}
&&{\rm E_{\alpha}}=\left\{x\in \R: \lim\limits_{t\longrightarrow0}
\left|U(t)f-f\right|>\alpha\right\}.\label{4.03}
\end{eqnarray}
Obviously, $E_{\alpha}\subset E_{1\alpha}\cup E_{2\alpha}$,
\begin{eqnarray}
&&\hspace{-1cm}{\rm E_{1\alpha}}=\left\{x\in \R: \lim\limits_{t\longrightarrow0}\left|U(t)
g-g\right|>\frac{\alpha}{2}\right\},{\rm E_{2\alpha}}=\left\{x\in \R: \lim\limits_{t\longrightarrow0}\left|U(t)h
-h\right|>\frac{\alpha}{2}\right\}.\label{4.05}
\end{eqnarray}
Obviously,
\begin{eqnarray}
E_{\alpha}\subset E_{1\alpha}\cup E_{2\alpha}.\label{4.06}
\end{eqnarray}
From Lemma 2.6,  we have
\begin{eqnarray}
\left|E_{1\alpha}\right|=0.\label{4.07}
\end{eqnarray}
Obviously,
\begin{eqnarray}
E_{2\alpha}\subset E_{21\alpha}\cup E_{22\alpha},\label{4.08}
\end{eqnarray}
where
\begin{eqnarray}
&&E_{21\alpha}=\left\{x\in \R: \sup\limits_{t>0}
\left|U(t)h\right|>\frac{\alpha}{4}\right\},E_{22\alpha}=\left\{x\in \R:
\left|h\right|>\frac{\alpha}{4}\right\}\label{4.010}.
\end{eqnarray}
Thus, from  Lemma 2.4,    we have
\begin{eqnarray}
&&\hspace{-1.6cm}\left| E_{21\alpha}\right|=\int_{E_{21\alpha}}dx\leq \int_{ E_{21\alpha}}
\frac{\left[\sup\limits_{t>0}\left|U(t)h\right|\right]^{p}}
{\alpha^{p}}dx\leq \frac{\left\|U(t)h\right\|_{L_{x}^{p}L_{t}^{\infty}}^{2r}}{\alpha^{p}}
\leq C\frac{\|h\|_{\hat{H}^{\frac{1}{p},\frac{p}{2}}}^{p}}{\alpha^{p}}
\leq C\frac{\epsilon^{p}}{\alpha^{p}}\label{4.011}.
\end{eqnarray}
By using the Hausdorff-Young inequality, we have
\begin{eqnarray}
&&\left| E_{22\alpha}\right|=\int_{E_{22\alpha}}dx
\leq \int_{ E_{22\alpha}}
\frac{\left|h\right|^{\frac{p}{2}}}
{\alpha^{\frac{p}{2}}}dx\leq \frac{\left\|h\right\|_{L_{x}^{\frac{p}{2}}}^{\frac{p}{2}}}{\alpha^{\frac{p}{2}}}
\leq C\frac{\|h\|_{\hat{L}^{\frac{p}{2}}}^{\frac{p}{2}}}{\alpha^{\frac{p}{2}}}
\leq C\frac{\epsilon^{\frac{p}{2}}}{\alpha^{\frac{p}{2}}}.\label{4.012}
\end{eqnarray}
From  (\ref{4.07}), (\ref{4.011}) and (\ref{4.012}),  we have
\begin{eqnarray}
&&\left|E_{\alpha}\right|\leq\left|E_{1\alpha}\right|+\left|E_{2\alpha}\right|\leq \left|E_{1\alpha}\right|+\left|E_{21\alpha}\right|+\left|E_{22\alpha}\right|\leq \frac{C\epsilon^{p}}{\alpha^{p}}+\frac{C\epsilon^{\frac{p}{2}}}{\alpha^{\frac{p}{2}}}\label{4.013}.
\end{eqnarray}
Thus, since  $\epsilon$  is arbitrary, from (\ref{4.013}),   we have
\begin{eqnarray}
\left|E_{\alpha}\right|=0.\label{4.014}
\end{eqnarray}
Thus, from  (\ref{4.03}) and  (\ref{4.014}), we have
$
U(t)f-f\longrightarrow 0(t\longrightarrow 0)\label{4.015}
$
almost everywhere.

This completes the proof of Theorem 1.1.

\bigskip

\noindent {\large\bf 5. Proofs of Theorem  1.2, 1.3}

\setcounter{equation}{0}

 \setcounter{Theorem}{0}

\setcounter{Lemma}{0}

\setcounter{section}{5}

By using (\ref{2.05})-(\ref{2.06}),  Lemmas 2.5, 2.6 and a proof similar to Theorem 1.1, we can obtain Theorems 1.2, 1.3.

\bigskip

\noindent {\large\bf 6. Proof of Theorem  1.4}

\setcounter{equation}{0}

 \setcounter{Theorem}{0}

\setcounter{Lemma}{0}

\setcounter{section}{7}

In this section, we  present the counterexamples showing that $s\geq \frac{1}{p}$
is the necessary condition for the maximal function  estimate.

 \noindent{\bf Proof of Theorem  1.4.}
 We choose $f=\frac{1}{\sqrt{2\pi}}\int_{\SR}e^{ix\xi}2^{-k(s+\frac{1}{(\frac{p}{2})^{\prime}})}\chi_{2^{k}\leq |\xi|\leq 2^{k+1}}(\xi)d\xi$, obviously,
 \begin{eqnarray}
 \|f\|_{\hat{H}^{s,\frac{p}{2}}}\sim 1.\label{6.01}
 \end{eqnarray}
We choose $t\leq \frac{\delta}{100}2^{-2k}$, where $\delta $  will be chosen later.
 For $|x|\leq 2^{-k}$ and sufficiently small $\delta$,
we have
\begin{eqnarray}
\|U(t)f\|_{L_{x}^{p}L_{t}^{\infty}}\sim 2^{-k(s-\frac{1}{2}+\frac{1}{p})}\label{6.02}.
\end{eqnarray}
From
$
\|U(t)f\|_{L_{x}^{p}L_{t}^{\infty}}
\leq C\|f\|_{\hat{H}^{s,\frac{p}{2}}(\SR)}\label{6.03}
$
and (\ref{6.01})-(\ref{6.02}), we have
\begin{eqnarray}
2^{-k(s-\frac{1}{p})}\leq C\label{6.04}.
\end{eqnarray}
We know that
 for sufficiently large $k,$ when $s<\frac{1}{p}$, (\ref{6.04}) is invalid.

This completes the proof of Theorem 1.4.
\bigskip

\noindent {\large\bf 7. Proof of Theorem  1.5}

\setcounter{equation}{0}

 \setcounter{Theorem}{0}

\setcounter{Lemma}{0}

\setcounter{section}{7}
In this section, we apply  Lemmas 3.2-3.4 to establish Theorem 1.5.

\noindent{\bf Proof of Theorem  1.5.}
We first show that if $f$ is a rapidly decreasing function,  $\forall\alpha>0$, such that
\begin{eqnarray}
&&\lim\limits_{t\rightarrow0}\mathbb{P}\left(\omega\in\Omega:\left|U(t)f^{\omega}(x)
-f^{\omega}(x)\right|>\alpha\right)=0.\label{7.01}
\end{eqnarray}
From Lemma 3.2,  $\forall\epsilon>0$ such  that $Ce\epsilon(\ln\frac{C_{1}}{\epsilon})^{\frac{1}{2}}<\alpha$, we have
\begin{eqnarray}
&&\mathbb{P}\left(\left\{\omega\in\Omega:|U(t)f^{\omega}-f^{\omega}|>\alpha\right\}\right)
\leq C_{1}e^{-\left(\frac{\alpha}{C|t|e}\right)^{2}}.\label{7.02}
\end{eqnarray}
It is easy to see that for $|t|<\epsilon$, we have
\begin{eqnarray}
&&\mathbb{P}\left(\left\{\omega\in\Omega:|U(t)f^{\omega}-f^{\omega}|>\alpha\right\}\right)
\leq C_{1}e^{-\left(\frac{\alpha}{Ce\epsilon}\right)^{2}}\leq\epsilon.\label{7.03}
\end{eqnarray}
Thus, from (\ref{7.03}), we derive  that (\ref{7.01}) is valid.
Next we consider $f\in\hat{L}^{r}(\R^{n}),2\leq r<\infty$. When $f\in\hat{L}^{r}(\R^{n}),2\leq r<\infty$,
  from  Lemma 2.5, we have that   $\forall\epsilon>0$,
 there exists a rapidly decreasing function $g$ such that $f=g+h$ which means
 $f^{\omega}=g^{\omega}+h^{\omega}$, where $\|h\|_{\hat{L}^{r}}<\epsilon$. Thus, we have
\begin{eqnarray}
&&\left|U(t)f^{\omega}-f^{\omega}\right|\leq\left|U(t)g^{\omega}-g^{\omega}\right|
+\left|U(t)h^{\omega}-h^{\omega}\right|.\label{7.04}
\end{eqnarray}
We define
\begin{eqnarray}
&&\Omega_{5}^{c}=\left\{\omega\in\Omega:|U(t)f^{\omega}-f^{\omega}|>\alpha\right\},\Omega_{6}^{c}=\left\{\omega\in\Omega:|U(t)g^{\omega}-g^{\omega}|>\frac{\alpha}{2}\right\},\label{7.05}\\
&&\Omega_{7}^{c}=\left\{\omega\in\Omega:|U(t)h^{\omega}-h^{\omega}|>\frac{\alpha}{2}\right\},\Omega_{8}^{c}=\left\{\omega\in\Omega:|U(t)h^{\omega}|>\frac{\alpha}{4}\right\},\label{7.08}\\
&&\Omega_{9}^{c}=\left\{\omega\in\Omega:|h^{\omega}|>\frac{\alpha}{4}\right\}.\label{7.09}
\end{eqnarray}
From  (\ref{7.04}) and  (\ref{7.05}), we have
\begin{eqnarray}
&&\Omega_{5}^{c}\subset\Omega_{6}^{c}\cup\Omega_{7}^{c}.\label{7.010}
\end{eqnarray}
From Lemma 3.2,  $\forall \epsilon>0$, when $|t|<\epsilon$, we have
\begin{eqnarray}
&&\mathbb{P}\left(\Omega_{6}^{c}\right)\leq C_{1}e^{-\left(\frac{\alpha}{C|t|e}\right)^{2}}\leq  C_{1}e^{-\left(\frac{\alpha}{Ce\epsilon}\right)^{2}}.\label{7.011}
\end{eqnarray}
Obviously,
\begin{eqnarray}
&&\Omega_{7}^{c}=\Omega_{8}^{c}\cup\Omega_{9}^{c}.\label{7.012}
\end{eqnarray}
From Lemmas 3.3,   3.4,  since $\|h\|_{\hat{L}^{r}}<\epsilon$, we have
\begin{eqnarray}
\mathbb{P}\left(\Omega_{8}^{c}\right)\leq C_{1}e^{-\left(\frac{\alpha}{Ce\|h\|_{\hat{L}^{r}}}\right)^{2}}\leq C_{1}e^{-\left(\frac{\alpha}{Ce\epsilon}\right)^{2}},
\mathbb{P}\left(\Omega_{9}^{c}\right)\leq C_{1}e^{-\left(\frac{\alpha}{Ce\|h\|_{\hat{L}^{r}}}\right)^{2}}\leq C_{1}e^{-\left(\frac{\alpha}{Ce\epsilon}\right)^{2}}.\label{7.013}
\end{eqnarray}
From (\ref{7.012}) and  (\ref{7.013}),  we can obtain
\begin{eqnarray}
&&\mathbb{P}\left(\Omega_{7}^{c}\right)\leq \mathbb{P}\left(\Omega_{8}^{c}\right)+\mathbb{P}\left(\Omega_{9}^{c}\right)\leq 2C_{1}e^{-\left(\frac{\alpha}{Ce\epsilon}\right)^{2}}.\label{7.015}
\end{eqnarray}
Using (\ref{7.010}), (\ref{7.011}), (\ref{7.015}),  $\forall\alpha>0$,  $\forall\epsilon>0$ such that $Ce\epsilon(\ln\frac{C_{1}}{\epsilon})^{\frac{1}{2}}<\alpha)$,
 then we have
\begin{eqnarray}
&&\hspace{-1.7cm}\mathbb{P}\left(\Omega_{5}^{c}\right)\leq \mathbb{P}\left(\Omega_{6}^{c}\right)+\mathbb{P}\left(\Omega_{7}^{c}\right)\leq C_{1}e^{-\left(\frac{\alpha}{Ce\epsilon}\right)^{2}}+2C_{1}e^{-\left(\frac{\alpha}
{Ce\epsilon}\right)^{2}}\leq 3C_{1}e^{-\left(\frac{\alpha}{Ce\epsilon}\right)^{2}}\leq 3\epsilon.\label{7.016}
\end{eqnarray}
From (\ref{7.016}), we have
$\lim\limits_{t\longrightarrow 0}\mathbb{P}\left(\left|U(t)f^{\omega}(x)-f^{\omega}(x)\right|>\alpha\right)=0.
$ The proof is completed.

\bigskip

\leftline{\large \bf Acknowledgments}
 This work is supported by NSFC grants (No. 11401180) and  the education department of Henan Province under
grant number 21A110014.


\bigskip

\baselineskip=18pt

\leftline{\large\bf  References}



\begin{thebibliography}{99}


\bibitem{BOP2015} A. B\'{e}nyi, T. Oh and  O. Pocovnicu, On the probabilistic Cauchy theory of the cubic
 nonlinear Schr\"odinger equation on $\R^d, d\geq3,$ {\it Trans.
Amer. Math. Soc. Ser. B} 2(2015), 1-50.

\bibitem{BOP}
A. B\'{e}nyi, T. Oh and  O. Pocovnicu,
Wiener randomization on unbounded domains and an application to
almost sure well-posedness of NLS,
in: Excursions in Harmonic Analysis,
Vol. 4, Birkh\"{a}user/Springer, Cham (2015), pp. 3-25.


\bibitem{Bourgain1994}J. Bourgain,
Periodic nonlinear Schr\"{o}dinger equation and invariant measures, {\it Comm. Math. Phys.} 166(1994),  1-26.



\bibitem{Bourgain1995}J. Bourgain,  Some new estimates on osillatory integrals,
  In: Essays on Fourier Analysis in Honor of
Elias M. Stein, Princeton, NJ 1991. Princeton Mathematical Series, vol. 42,
 pp. 83.112. Princeton
University Press, New Jersey (1995).

\bibitem{Bourgain1996}J. Bourgain,
Invariant measures for the $2D$-defocusing nonlinear Schr\"{o}dinger equation, {\it  Comm. Math. Phys.} 176(1996),  421-445.







\bibitem{Bourgain2013}J. Bourgain,  On the Schr\"odinger maximal
 function in higher dimensions,  {\it Proc. Steklov Inst. Math.} 280(2013), 46-60.

\bibitem{Bourgain2016} J. Bourgain,  A note on the Schr\"odinger
maximal function,  {\it J. Anal. Math.} 130(2016),  393-396.


\bibitem{BTL}N. Burq and  N.  Tzvetkov,  Random data Cauchy theory for
 supercritical wave equations,  I. Local theory, {\it Invent.Math.} 173(2008),  449-475.

\bibitem{BTG}N. Burq and  N.  Tzvetkov, Random data Cauchy theory for supercritical wave equations,
 II. A global existence result, {\it Invent.Math.} 173(2008),  477-496.






\bibitem{CVV}  T. Cazenave,L.  Vega, M. C. Vilela,
A note on the nonlinear Schr\"odinger equation in weak $L^{p}$
spaces, {\it Commun. Contemp. Math.} 3(2001),  153-162.

\bibitem{Carleson} L. Carleson,  Some analytical problems related to statistical
 mechanics. Euclidean Harmonic Analysisi.
Lecture Notes in Mathematics, vol. 779, pp. 5.45, Springer, Berlin, (1979).

\bibitem{CG}Y. Chen, H. Gao,
 The Cauchy problem for the Hartree equations under random influences, {\it J. Diff. Eqns.} 259(2015), 5192-5219.



\bibitem{CLS}E. Compaan, R. Luc\'{a}, G. Staffilani ,
Pointwise Convergence of the Schr\"{o}dinger Flow, arXiv:1907.11192v1 [math.AP] 25 Jul 2019,  doi: 10.1093/imrn/rnaa036..





\bibitem{CLV}C. Cho, S. Lee and  A. Vargas, Problems on pointwise convergence of solutions
to the Schr\"odinger
equation,  {\it J. Fourier Anal. Appl.} 18(2012), 972-994.









\bibitem{Cowling}M.  Cowling,  Pointwise behavior of solutions to
 Schr\"odinger equations. In: Harmonic Analysis (Cortona,
1982). Lecture Notes in Mathematics, vol. 992, pp. 83.90. Springer, Berlin, (1983).

 \bibitem{DK} B. E. Dahlberg and  C. E. Kenig, A note on the almost everywhere behavior
  of solutions to the Schr\"odinger
equation. In: Proceedings of Italo-American Symposium in Harmonic Analysis,
University of Minnesota.
Lecture Notes in Mathematics, vol. 908, pp. 205.208. Springer, Berlin, (1982).

\bibitem {DG} C. Demeter and  S. Guo,  Schr\"odinger maximal function estimates via
 the pseudoconformal transformation,
 arXiv: 1608.07640.

\bibitem{DLM}B. Dodson, J. L\"{u}hrmann, D. Mendelson,
 Almost sure local well-posedness and scattering for the 4D cubic nonlinear Schr\"odinger equation, {\it Adv. Math.} 347(2019), 619-676.


\bibitem{Du} X.  Du, A sharp  Schr\"odinger  maximal estimate in $\R^{2}$,  {\it Dissertation,} 2017.

\bibitem{DGL} X. Du, L. Guth and   X.  Li,  A sharp Schr\"odinger maximal estimate in $\R^2,$
  {\it Ann. Math.} 188(2017), 607-640.


\bibitem{DZ}X.  Du and  R.  Zhang, Sharp $L^2$ estimates of the Schr\"odinger maximal function
 in higher dimensions, {\it Ann. Math.} 189(2019),
837-861.

\bibitem{DGZ}X.  Du, L. Guth,  X. Li and R. Zhang, Pointwise convergence of Schr\"{o}dinger solutions
 and multilinear refined Strichartz estimates, {\it Forum Math.Sigma} 6(2018).


 \bibitem{EL} D. Eceizabarrena and R. Luc$\grave{a}$, Convergence over fractals
 for the periodic Schr\"odinger equation,
arXiv:2005.07581.





\bibitem{F}  C.  Fefferman,  Inequalities for strongly singular convolution operators, {\it  Acta Math.} 124(1970), 9-36.


\bibitem{GS}G. Gigante and  F. Soria,  On the the boundedness in $H^{1/4}$ of the maximal square
 function associated with the Schr\"odinger equation,  {\it J. Lond. Math. Soc.} 77(2008), 51-68.







\bibitem{G}A.  Gr\"unrock, An improved local well-posedness result for the modified KdV equation,
{\it  Int. Math. Res. Not.}  61(2004),  3287-3308.








\bibitem{KPV1991}  C. E. Kenig,  G. Ponce and  L. Vega, Oscillatory integrals and
 regularity of dispersive equations,
{\it  India. Univ. Math. J.} 40(1991),  33-69.

\bibitem{KMV}R. Killip, J. Murphy, M. Visan,
Almost sure scattering for the energy-critical NLS with radial data below $H^{1}(\R^{4})$, {\it  Commun. Partial Diff. Eqns.} 44(2019), 51-71.


\bibitem{KMV2020}R. Killip, J.  Murphy, M. Visan, Invariance of white noise for KdV on the line, {\it  Invent. Math.} 222(2020),   203-282.


\bibitem{LRS}J. Lebowitz, H.  Rose, E. Speer,
 Statistical mechanics of the nonlinear Schr\"{o}dinger equation, {\it J. Statist. Phys.} 50(1988),  657-687.



\bibitem{Lee}S.  Lee,  On pointwise convergence of the solutions to
 Schr\"odinger equation
 in $\R^2,$ {\it Int. Math. Res. Not.}    (2006), 32597.







\bibitem{LR2015}R. Luca and  M. Rogers,  An improved neccessary condition for Schr\"odinger
 maximal estimate,  arXiv:1506.05325.

\bibitem{LR2017}R. Luca and  M. Rogers,  Coherence on fractals versus pointwise convergence
 for the Schr\"odinger equation,
{\it Commun. Math. Phys.} 351(2017), 341-359.

\bibitem{LMCPDE}J. L\"uhrmann and  D. Mendelson,  Random data Cauchy
 theory for nonlinear wave equations of power-type on
$\R^3$, {\it Comm. Partial Diff. Eqns.}  39(2014),  2262-2283.


\bibitem{MYZ2015} C. Miao, J. Yang and  J. Zheng, An improved maximal inequality for 2D
 fractional order Schr\"odinger
operators,  Stud. Math. 230(2015), 121-165.

\bibitem{MZZ2015}C.  Miao, J. Zhang and  J. Zheng, Maximal estimates for Schr\"odinger
 equation with inverse-square potential,
{\it Pac. J. Math.} 273(2015), 1-19,

\bibitem{M}J. Murphy,
Random data final-state problem for the mass-subcritical NLS in $L^{2}$, {\it Proc. Amer. Math. Soc.} 147(2019), 339-350.





\bibitem{MVV}A.  Moyua, A. Vargas and  L. Vega,  Schr\"odinger maximal function and
restriction properties of the Fourier
transform,  {\it Int. Math. Res. Not.} 1996(1996), 793-815.

\bibitem{MV} A. Moyua and L. Vega,
Bounds for the maximal function associated to periodic solutions of one-dimensional
 dispersive equations,
{\it Bull.  Lon. Math. Soc.}  40(2008),  117-128.



\bibitem{OOP}T. Oh, M. Okamoto, O. Pocovnicu,
 On the probabilistic well-posedness of the nonlinear Schr\"{o}dinger equations with non-algebraic nonlinearities, {\it  Discrete Contin. Dyn. Syst.} 39(2019),  3479-3520.

\bibitem{OP}T. Oh, O. Pocovnicu,
Probabilistic global well-posedness of the energy-critical defocusing quintic nonlinear wave equation on $\R^{3}$,
J. Math. Pures Appl. 105(2016),  342-366.










\bibitem{Shao} S. Shao, On localization of the Schr\"odinger maximal operator,
 arXiv: 1006.2787v1.

\bibitem{S} P. Sj\"olin,  Regularity of solutions to the Schr\"odinger equation,
{\it  Duke Math. J.}  55(1987), 699-715.


\bibitem {Tao} T. Tao,  A sharp bilinear restriction estimate for parabloids,
{\it  Geom. Funct. Anal.} 13(2003), 1359-1384.

\bibitem{TV}T. Tao and  A. Vargas, A bilinear approach to cone multipliers,
II. Appl.{\it  Geom. Funct. Anal.} 10(2003), 216-258.



\bibitem{Vega}L. Vega,  Schr\"odinger equations: pointwise convergence to
the initial data,  {\it Proc. Am. Math. Soc.}
102(1988), 874-878.





\bibitem{WZ}X. Wang and  C. J. Zhang,  Pointwise convergence of solutions to the
 Schr\"odinger equation on manifolds,
{\it  Canad. J. Math. }  71(2019),  983-995.





\bibitem{Zhang}C. Zhang,   Pointwise convergence of solutions to Schr\"odinger type equations,
{\it Nonli. Anal. TMA}109(2014), 180-186.



\bibitem{ZF2011}T. Zhang and   D. Fang,  Random data Cauchy theory for the incompressible three dimensional Navier-Stokes equations,
{\it Proc. Amer. Math. Soc.} 139(2011), 2827-2837.

\bibitem{ZF2012}T. Zhang and   D. Fang,  Random data Cauchy theory for the generalized incompressible Navier-Stokes equations,
{\it J. Math. Fluid Mech.} 14(2012),  311-324.




\end{thebibliography}
\end{document}